\documentclass[12pt]{article}
\usepackage{amssymb}
\usepackage{multirow}

\def\a{\alpha}
\def\b{\beta}

\def\e{\varepsilon}

\def\g{\gamma}

\def\i{\iota}

\def\l{\lambda}
\def\o{\omega}

\def\t{\tau}
\def\u{\upsilon}

\def\G{\Gamma}

\def\O{\Omega}

\chardef\tempcat=\the\catcode`\@
\catcode`\@=11
\def\cyracc{\def\u##1{\if \i##1\accent"24 i
    \else \accent"24 ##1\fi }}
\newfam\cyrfam

\font\cmsslll=cmss10 at 14 pt

\DeclareFontFamily{OT1}{msb}{}{}
\DeclareFontShape{OT1}{msb}{m}{n}
 {  <5> <6> <7> <8> <9> <10> gen * msbm
      <10.95><12><14.4><17.28><20.74><24.88>msbm10}{}
\DeclareMathAlphabet{\bubble}{OT1}{msb}{m}{n}

\def\bO{{\mathbb O}}

\def\bR{{\mathbb R}}

\def\bC{{\mathbb C}}
\def\bH{{\mathbb H}}

\newfont{\goth}{eufm10 scaled \magstep1}
\def\ga{\mbox{\goth a}}

\def\gf{\mbox{\goth f}}
\def\gg{\mbox{\goth g}}
\def\gh{{\mbox{\goth h}}}

\def\gl{\mbox{\goth l}}
\def\gm{\mbox{\goth m}}

\def\gp{\mbox{\goth p}}

\def\gu{\mbox{\goth u}}

\def\gso{\mbox{\goth so}}
\def\gsu{\mbox{\goth su}}

\def\gsp{\mbox{\goth sp}}

\def\ggl{\mbox{\goth gl}}

\newfont{\mcal}{eusm10 scaled \magstep1}

\def\cx{\mbox{\mcal X}}

\def\p{\partial}

\def\ra{\rightarrow}
\def\Id{\mathrm{Id\;}}

\def\Ad{\mathrm{Ad\;}}
\def\ad{\mathrm{ad\;}}

\newtheorem{Th}{Theorem}
\newtheorem{MTh}{Main Theorem}
\newtheorem{Prop}{Proposition}
\newtheorem{Cor}{Corollary}
\newtheorem{Lem}{Lemma}
\newtheorem{Def}{Definition}

\def\bt{\begin{Th}}
\def\et{\end{Th}}
\def\bp{\begin{Prop}}
\def\ep{\end{Prop}}
\def\bc{\begin{Cor}}
\def\ec{\end{Cor}}
\def\bl{\begin{Lem}}
\def\el{\end{Lem}}
\def\bd{\begin{Def}}
\def\ed{\end{Def}}
\def\pf{\noindent{\it Proof. }}
\def\qed{\hspace{2ex} \hfill $\square $ \par \medskip}

\def\be{\begin{equation}}
\def\ee{\end{equation}}
\def\arr{\begin{array}{rlll}}
\def\ea{\end{array}}
\def\bea{\begin{eqnarray}}
\def\eea{\end{eqnarray}}
\def\bean{\begin{eqnarray*}}
\def\eean{\end{eqnarray*}}

\catcode`@=11
\@addtoreset{equation}{section}
\catcode`@=12


\begin{document}
\begin{center}
{\LARGE On the Geometry of Spaces of Oriented Geodesics}\\
\vskip 1.0 true cm {\cmsslll Dmitri V.\ Alekseevsky$\,$,}\\
{\small School of Mathematics and Maxwell Insitute for
Mathematical Sciences, The Kings Buildings, JCMB, University of
Edinburgh,
Mayfield  Road, Edinburgh, EH9 3JZ,UK. \\
D.Aleksee@ed.ac.uk}
\vskip 1.0 true cm {\cmsslll Brendan Guilfoyle$\,$,}\\
{\small Department of Computing and Mathematics,
IT Tralee, Clash, Tralee,
County Kerry, Ireland.  \\
brendan.guilfoyle@ittralee.ie}
\vskip 1.0 true cm {\cmsslll Wilhelm Klingenberg$\,$.}\\
{\small School of Mathematical Sciences, University of Durham, Durham DH1
3LE, UK. \\
wilhelm.klingenberg@dur.ac.uk}

\vspace{0.5in}

{\bf Abstract}:
Let $M$ be either a  simply connected  pseudo-Riemannian space of constant curvature or a rank one Riemannian symmetric space other than $\mathbb{O}H^2$, and consider the space
$L(M)$ of oriented geodesics of $M$. The  space $L(M)$  is a smooth
homogeneous manifold  and in this paper we describe all invariant
 symplectic structures, (para)complex structures, pseudo-Riemannian metrics   and
(para)K\"ahler  structure on $L(M)$.
\end{center}
{\small

 }

\vfill \hrule width 3.cm \vskip 0.2 true cm {\small \noindent This
work was supported by Leverhulme Trust, EM/9/2005/0069.}
\newpage
\tableofcontents

\section{Introduction}
\subsection{Background}

The  geometry  of the set of  straight lines  of projective space $P^3$ and Euclidean  space $E^3$ is a classical
  subject of investigations of such  19th century geometers  as Grassmann,  Pl\"ucker, F. Klein and Study.
They studied natural correspondences between submanifolds (i.e. points, curves and surfaces) in $E^3$   and  submanifolds of $L(E^3)$, the
space of oriented lines of $E^3$.
For example, a  point  $ p \in E^3$ defines a  surface in $L(E^3)$  which  consists of the oriented lines through $p$, a curve $C\subset E^3 $
 defines  three curves in $L(E^3)$  associated with  the  Frenet frame
of $C$, and a  surface $S \subset E^3$  defines a surface in $L(E^3)$ by its oriented normal lines. Conversely, distinguished  (from the point of  view of
intrinsic geometry of $L(E^3)$)  submanifolds  determine  special  families of oriented lines in $E^3$.

 In \cite{study}, Study  identified  the  space $L(E^3)$ of oriented lines in $E^3$ with the  2-sphere (the Study sphere) over the dual numbers
 and defined and studied  a  notion of   distance between oriented lines.
 For a nice modern exposition  and generalization of these results with applications to computational geometry, computer graphics
 and visualization, see \cite{PW}.

 A natural complex  structure in the space $L(E^3)$ of oriented lines  has been considered by Hitchin, who used it for a description of  monopoles \cite{hitch}.
 In \cite{GK1}, two of the authors defined  a neutral K\"ahler  structure in the space $L(E^3)$ and gave its geometric description.
Recently, it has been used in the solution of a long-standing conjecture of Carath\'eodory \cite{GK2}.

In general, however, the space $L(M)$ of oriented geodesics of a (complete) Riemannian manifold  is not  a smooth manifold  and has very bad topology.
But $L(M)$ is a smooth manifold if $(M,g)$ is either an Hadamard manifold (i.e. complete simply connected  Riemannian manifold of non-positive curvature)
or a manifold with closed geodesics of the same length ({\it aufwiederseen} manifolds).
Symplectic and Riemannian structures in the space  of geodesics $L(M)$ of aufwiederseen manifolds  are discussed by Besse
\cite{B}, see also  Reznikov \cite{R} \cite{R1}.
The symplectic form on the space of geodesics of a Hadamard manifold  is described  in \cite{ferr}.

Geometric  structures in the space of oriented geodesics of hyperbolic 3-space are studied in \cite{GG}, while  Salvai has addressed the existence
and uniqueness of pseudo-Riemannian metrics in the spaces $L(E^n)$  and $L(H^n)$ \cite{S2} \cite{S1}. Note that both spaces are  homogeneous  manifolds
of the corresponding isometry group.  Salvai  proved  that $L(E^n)$  admits a pseudo-Riemannian metric invariant under
a transitive subgroup $G$  of the isometry group $I(E^n) = E(n) = SO(n) \cdot E^n$  only for $n=3,7$ and gave  an explicit description of the
corresponding metric.

\subsection{Main results}

The aim of the paper is to  describe  the natural  geometric structures
in the  space $L(M)$ of oriented geodesics of $M$,  where $M$ is either a simply connected  pseudo-Riemannian
space of constant curvature,   or a  rank one Riemannian   symmetric space other than $\mathbb{O}H^2$. In these cases  the  space of geodesics $L(M)$  is a smooth
homogeneous manifold  and  we use Lie groups and Lie algebras  to describe all invariant
 symplectic structures, (para)complex structures, pseudo-Riemannian metrics   and
(para)K\"ahler  structure on $L(M)$.

More specifically, let $S^{p,q}= \{ x \in E^{p+1,q}, x^2 = 1 \}$, where $E^{p+1,q}$ is ${\mathbb{R}}^{p+q+1}$ endowed with the flat metric of signature
$(p+1,q)$. Here, for $p = 0$ we
assume in addition that $x^0 > 0$  so that $S^{p,q}$ is connected. The induced metric on $S^{p,q}$
has signature $(p,q)$ and is of constant curvature $1$.

Let $L^+(S^{p,q})$ (respectively $L^-(S^{p,q})$)
be the set of spacelike (respectively, timelike) geodesics in $S^{p,q}$ and similarly for  $L^\pm(E^{p+1,q})$. We prove:

\begin{MTh} \label{MT:1}
For the flat pseudo-Euclidean spaces $E^{p+1,q}$
\begin{enumerate}
\item[i)]
The space $L^-(E^{p+1,q}) = E(p+1,q)/ SO(p,q) \cdot {\mathbb{R}^+}$ is  a symplectic symmetric space with an invariant
Grassmann structure defined by  a decomposition $\gm = W \otimes {\mathbb{R}}^2$. Moreover, if  $n = p+1 +q =3$,
it has an invariant K\"ahler structure $(g,J)$ of  neutral signature $(2,2)$. In addition, $L^+(E^{p+1,q})=L^-(E^{q,p+1})$.
\end{enumerate}
while for the non-flat constant curvature manifolds $S^{p,q}$
\begin{enumerate}
\item[ii)] Suppose that $ p+q > 3$. Then there exists a unique (up to  scaling)
invariant symplectic structure $\omega$  and a unique (up to a sign)
invariant complex structure $I^+ =J$ on $L^+(S^{p,q})$
and a unique (up to sign) invariant para-complex structure $K=I^-$  on  $L^-(S^{p,q})$.
  There exists unique (up to scaling) invariant
pseudo-Riemannian metric $g^\e = \o \circ I^\e$  on $L^\e(S^{p,q})$ which is K\"ahler of signature
 $(2(p-1),2q)$ for
$\e =+$ and para-K\"ahler (of neutral signature)  for $\e =-$.
\item[ iii)] Suppose that $p+q =3$.  Then there  are 2-linearly independent invariant (parallel and closed)
 2-forms $\o, \o'$ on $L^\e(S^{p,q})$ with values $\o_{\gm} = \o_H \otimes g_V, \, \o'_{\gm} = g_H \otimes \o_V$.
Any invariant metric has the form $h = \lambda g + \mu g'$ where $g'$ is a neutral metric with
value $g'_{\gm } = \o_H \otimes \o_V$.    Any metric $h$ is K\"ahler for $\e=1$ (respectively,
 para-K\"ahler for $\e =-$) with respect to the complex (respectively, para-complex)  structure
$I^\e = I^\e_H \otimes 1 $ with the K\"ahler form $h \circ I^\pm = \lambda \o + \mu \o'$.
Moreover, the endomorphism $I' = 1 \otimes I_V$ of $\gm$ defines  an  invariant parallel
 $h$-skew-symmetric  complex  structure of $L^\e(S^{p,q})$ if $\e =1, (p-1,q)= (2,0)$ or $(0,2)$ or
$\e =-$  and $(p,q-1) = (2,0)$ or $(0,2)$  and skew-symmetric parallel  para-complex structure
if $\e =+, (p-1,q)= (1,1)$ or $\e =-, (p,q-1) = (1,1)$. The K\"ahler or para-K\"ahler structure
$(h, I')$ has  the K\"ahler form $h \circ I' = \lambda \o + \mu \o'$.
\end{enumerate}

\end{MTh}

Consider now the rank one Riemannian symmetric spaces of non-constant
curvature. That is, $M = G/K$ is one of  the projective spaces
\[
\bC P^n =SU_{n+1}/ U_n, \qquad \bH P^n = Sp_{n+1}/ Sp_1 \cdot Sp_n,\qquad \bO P^2 = F_4/Spin_9
\]
or one of  the dual  hyperbolic spaces
\[
\bC H^n =SU_{1,n}/ U_n,\qquad \bH H^n = Sp_{1,n}/ Sp_1 \cdot Sp_n,\qquad \bO H^2 = F^{non-comp}_4/Spin_9 .
\]

For the projective spaces we prove:

\begin{MTh}\label{MT:2}
\begin{enumerate}

\item[i)]
The space $L(\bC P^n) =S U_{n+1}/T^2\cdot SU_{n-1}$  has a one--parameter
family $\omega^t = \omega_1 + t \omega_0$ of invariant symplectic forms (up to scaling) and four invariant
almost complex structures, up to sign, two of them being integrable.
All (almost) complex structures $J$ are compatible  with  $\omega^t$ i.e.  they define  an (almost)
K\"ahler or pseudo-K\"ahler metric $g=\omega^t\circ J$.
\item[ii)]
The spaces $L(\bH P^n) = Sp_{n+1}/T^1\cdot Sp_1\cdot Sp_{n-1}$ and  $L(\bO P^2)= F_4/T^1\cdot Spin_7$
has a unique (up to scaling) invariant symplectic forms $\omega$ and
unique (up to sign) invariant complex structure $J$.
The pair $(\omega, J)$ defines   a unique (up to scaling)  invariant K\"ahler metric
$g=\omega\circ J$.

\item[iii)] None of the above spaces have any invariant almost para-complex structures.
\item[iv)] The canonical
symplectic structure on $L(\bC P^n)$ is identified with $\o_1$.
\end{enumerate}
\end{MTh}

while for the hyperbolic spaces (other than $\mathbb{O}H^2$) we show

\begin{MTh}\label{MT:3}

\begin{enumerate}
\item[i)]
The space $L(\bC H^n)= SU_{1,n}/T^2\cdot SU_{n-1}$  has a one--parameter family of
invariant  symplectic structures $\omega^t = \omega_1 + t \omega_0$,
and two (up to sign) invariant almost para-complex structures $K^{\pm}$, one of them
being integrable; and both are consistent
with  $\omega^t$
i.e. $(K^\pm, \omega^t)$  defines a para-K\"ahler metric $g = \omega \circ K^\pm $.
\item[ii)]
The space $L(\bH H^n) = Sp_{1,n}/T^1\cdot Sp_1\cdot Sp_{n-1}$  admits a unique (up to scaling)
invariant symplectic form $\omega$ and  two (up to sign) invariant almost para-complex
structures, one of which is integrable. Both are compatible with $\omega$.
\item[iii)]
None of the above spaces have any invariant almost complex structures.
 \item[iv)]The canonical
symplectic structure on the space geodesics $L(\bC H^n)$ is $\o_1$.

\end{enumerate}
\end{MTh}

In Table 1 we summarize the results of Main Theorems \ref{MT:1}, \ref{MT:2} and \ref{MT:3}.
\newpage

\begin{table}
\caption{Invariant Geometric Structures}
\begin{tabular}{|c|c||c|c|c|c|c|c|c|c|c|}
\hline
\multicolumn{2}{|c|}{}  & Symplectic&  \multicolumn{2}{|c|}{Complex} & \multicolumn{2}{|c|}{Para-Complex} & \multicolumn{2}{|c|}{K\"ahler}& \multicolumn{2}{|c|}{ParaK\"ahler}   \\
\multicolumn{2}{|c|}{}&Structure & Int & Non & Int & Non & Int & Non & Int & Non \\\hline\hline
\multirow{2}{*}{$L^{\pm}(E^{p+1,q})$} &${\scriptstyle{p+1+q\neq 3}}$& \multirow{2}{*}{Symmetric}&$\emptyset$ &$\emptyset$ &$\emptyset$ &$\emptyset$ & $\emptyset$&$\emptyset$ &$\emptyset$ & $\emptyset$\\
   &${\scriptstyle{p+1+q=3}}$& & ${\mathbb R}$&$\emptyset$ &$\emptyset$ &$\emptyset$ &${\mathbb R}$ &$\emptyset$ &$\emptyset$ &$\emptyset$  \\\hline
\multirow{4}{*}{$L^+(S^{p,q})$}&${\scriptstyle{p+q>3}}$& 1&1 &$\emptyset$ & $\emptyset$ &$\emptyset$ &1 &$\emptyset$ &$\emptyset$ &$\emptyset$ \\
   &${\scriptstyle{(p,q)=(3,0)}}$& ${\mathbb R}$ &1 & $\emptyset$&$\emptyset$ &$\emptyset$ &${\mathbb R}$ & $\emptyset$& $\emptyset$& $\emptyset$\\
&${\scriptstyle{(p,q)=(1,2)}}$& ${\mathbb R}$ &1 & $\emptyset$&$\emptyset$&$\emptyset$ &${\mathbb R}$ & $\emptyset$& $\emptyset$& $\emptyset$\\
&${\scriptstyle{(p,q)=(2,1)}}$& ${\mathbb R}$ &$\emptyset$ & $\emptyset$&1 &$\emptyset$ &$\emptyset$ & $\emptyset$& ${\mathbb R}$& $\emptyset$\\\hline
\multirow{4}{*}{$L^-(S^{p,q})$}&${\scriptstyle{p+q>3}}$& 1 &$\emptyset$ &$\emptyset$ &1 & $\emptyset$& $\emptyset$& $\emptyset$& 1&$\emptyset$ \\
   &${\scriptstyle{(p,q)=(0,3)}}$& ${\mathbb R}$& 1&$\emptyset$ &$\emptyset$ &$\emptyset$ &${\mathbb R}$ &$\emptyset$ &$\emptyset$ & $\emptyset$\\
   &${\scriptstyle{(p,q)=(2,1)}}$& ${\mathbb R}$& 1&$\emptyset$ &$\emptyset$ &$\emptyset$ &${\mathbb R}$ &$\emptyset$ & $\emptyset$& $\emptyset$\\
   &${\scriptstyle{(p,q)=(1,2)}}$& ${\mathbb R}$& $\emptyset$&$\emptyset$ & 1&$\emptyset$ &$\emptyset$ &$\emptyset$ & ${\mathbb R}$& $\emptyset$\\\hline
\multicolumn{2}{|c||}{$L({\mathbb C}P^n)$}& ${\mathbb R}$ & 2& 2 &$\emptyset$ &$\emptyset$ &2 &2 &$\emptyset$ &$\emptyset$ \\\hline
\multicolumn{2}{|c||}{$L({\mathbb H}P^n)$}& 1& 1 &$\emptyset$ &$\emptyset$ & $\emptyset$&1 &$\emptyset$ &$\emptyset$ &$\emptyset$ \\\hline
\multicolumn{2}{|c||}{$L({\mathbb O}P^2)$}& 1& 1 &$\emptyset$ &$\emptyset$ & $\emptyset$&1 &$\emptyset$ &$\emptyset$ &$\emptyset$ \\\hline
\multicolumn{2}{|c||}{$L({\mathbb C}H^n)$}& ${\mathbb R}$&$\emptyset$ &$\emptyset$ & 1& 1&$\emptyset$ &$\emptyset$ & 1&1 \\\hline
\multicolumn{2}{|c||}{$L({\mathbb H}H^n)$}& 1            &$\emptyset$ &$\emptyset$ & 1& 1&$\emptyset$ &$\emptyset$ &1 &1 \\\hline
\end{tabular}
\end{table}

\subsection{Outline of paper}

This paper is organised as follows. In the following section we consider the set of oriented geodesics $L= L(M)$  of a general complete
Riemannian manifold $M$.
 Though  $L$ may not be  a manifold,  we can still define
differential geometric objects on $L$ in terms of $G_1$-invariant objects on $SM$, where $SM$ is the unit sphere bundle and $G_1= \{exp t \G\}$ is the
geodesic flow. For example, the algebra of  smooth functions ${\mathcal F}(L)$ on $L$ is defined  as the algebra
of $G_1$-invariant functions on $SM$ (that is  first integrals of the geodesic flow).
We  sketch  this  approach in section 2 and define  the canonical  symplectic form on $L$. In the case when
$L$ is a manifold  this  coincides  with
the  standard  symplectic form defined  in  \cite{B}.

Section 3 considers the oriented geodesics of pseudo-Riemannian manifolds of constant curvature, taking the flat and non-flat cases separately.
Main Theorem \ref{MT:1} follows from Theorems \ref{t:mt1a} and \ref{t:mt1b} of sections 3.1 and 3.2, respectively.

In section 4 we turn to rank one Riemannian symmetric spaces with non-constant curvature, dealing separately with the complex projective
and hyperbolic spaces,
the quaternionic projective and hyperbolic spaces and the Cayley plane. In particular, Main Theorem \ref{MT:2}
follows from Theorem \ref{t:mt2a} of section 4.1, Theorem \ref{t:mt2b} of section 4.2 and Theorem \ref{t:mt2c} of section 4.3, while
Main Theorem \ref{MT:3}
follows from Theorem \ref{t:mt3a} of section 4.1 and Theorem \ref{t:mt3b} of section 4.2.


\section{The Space of Oriented Geodesics on a Riemannian Manifold}

\subsection{Tangent and unit sphere bundle}

Let $M$ be an $n$-dimensional manifold and $\pi :TM \to M $ be the tangent bundle of $M$. Local coordinates $(x^i)$ on $M$ give rise
to local coordinates
$(x^i,y^i)$ on $TM$ via $v_x = y^i\p/\p x^i$. Then $T_{(x,y)}TM \ni (x,y,\dot x, \dot y)$.
The vertical subspace $V_{(x,y)}TM$  is given by $(x,y,0,\dot y)$ and we denote by $N$
the canonical endomorphism
\[
N : T_{(x,y)}TM \to V_{(x,y)}TM \to  0 : \;\;
 (x,y,\dot x, \dot y) \to (x, y, 0,\dot x) \to 0,
\]
with $N^2 = 0$. Note that
\[  \mathrm{Ker} \; N = \mathrm{Im}  \; N = T^v(TM) = \{(x,y, 0,\dot y) \},
\]
is the vertical subbundle of $T(TM)$. The restriction of the projection $\pi_* : T(TM) \to TM$
onto the vertical subspace  is  an isomorphism
\[
  \pi_* : T^v_{(x,y)}(TM) \to T_xM,\,\, (x,y,\dot x, \dot y) \mapsto (x,\dot x).
\]
Assume now that $(M,g)$ is
Riemannian and $SM$ is the unit sphere bundle:
\[
SM = \{ (x,y) \in TM, \, g_{ij}(x)y^iy^j = 1 \}.
\]
The metric $g$ induces an isomorphism of the tangent bundle
onto cotangent bundle $T^*M  = \{\alpha = (x,p),\, \alpha = p_idx^i \}$ given by
\[
(x^i,y^i) \to (x^i,p_i = g_{ij}y^j) .
\]
The pull-back of the canonical
one-form $\a = p_i dx^i$ and symplectic form $\o = d \a$ of $T^*M$ are given by
$$\a^g = g^* \a = g_{ij}y^i dx^j,\,\,
 \o^g = g \circ \o = g\circ d \a = d(g_{ij}y^i)\wedge d x^j.$$
 The tangent  bundle of the unit sphere bundle has the canonical  decomposition
$T_{(x,y)}SM =  V_{(x,y)}SM + H^g_{(x,y)} SM$
 into vertical space and horizontal subspaces.

\bl
The sphere bundle of the Riemannian manifold $(M,g)$
admits the canonical contact structure
\[
 \theta = \a^g|SM,  \;\;\; \a^g = g^* \a = g_{ij}y^i dx^j,
\]
where the associated Reeb vector field is the geodesic vector field given by
$$\Gamma =  y^i\p/\p x^i - \Gamma^i_{jk}y^j y^k \p/\p y^i.$$
The horizontal lift of $ \p/\p x^i = \p_i $  into $SM$   is given by

 $$\nabla_i = \p_i  - \Gamma^k_{ij}(x)y^j\p/\p y^k, (x,y) \in SM.$$
In particular, $\Gamma$ preserves  $\theta$  and $ d \theta$.

\el

\pf
The last claims are verified as follows. A vector field along $x(t)$ given by $X(t) = y^i(x(t))\p_i $ is parallel if
$$
\dot{X}^i + \G^i_{jk}(x(t))\dot{x}^iX^j =0.
$$
Then $(x(t),X(t))$ is the horizontal lift of the curve $x(t)$. The horizontal space is spanned by
the tangent vectors of such lifts, namely
\[
H_{(x,y)} = {\rm{span}}  \{(\dot{x}^i, -\G^i_{jk}\dot{x}^jy^k) \}.
\]
The vector field $\Gamma =  y^i\p/\p x^i - \Gamma^i_{jk}y^j y^k \p/\p y^i$ on  $TM$ is tangent to $SM$ since
$$
\nabla_k (g_{ij}y^iy^j) = g_{ij,k} - g_{im}\G^m_{kl}  - g_{jm}\G^m_{ki} = g_{ij,k} - \G_{i,kj}  - \G_{j,ki} =0.
$$
Moreover, $\G$ is a geodesic vector field since its integral curves satisfy the geodesic equation. Therefore we have
the decomposition into vertical and horizontal parts
$$
T_{(x,y)}SM = V_{(x,y)}SM + H_{(x,y)}SM = \{   v^i\p/\p y^i : g_{ij}y^i v^j  = 0 \} + \{u^k \nabla_k  \}.
$$
We compute
\[
d\theta = g_{ij,k}y^j dx^i\wedge dx^k + g_{ij}dy^j \wedge dx^i.
\]
One easily checks that $\theta(\G) = 1$ and $ d\theta (\G)=0 $.
Clearly the form $d( g^* \a) = g^* d\a = g^* \o $ is non-degenerate on $TM$. Therefore its restriction
to $SM$ has one-dimensional kernel spanned by $\G$ and so $\G$ is the Reeb vector field of the
contact form $\theta$  and it preserves $\theta$ and $d \theta$.
\qed

We associate with the vector field $X = X^i \p_i$ the function $f_X$ on $SM$ given by
$f_X(x,y)= g_{ij}X^i y^j$.
\bl
The covariant derivative $\nabla_i X$ corresponds to the Lie derivative of $f_X$ in direction of the vector field
$\nabla_i$ :
$$
f_{\nabla_i X}= \nabla_i f_X.
 $$
\el
\pf
Applying the vector field $\nabla_i$ to the function $f_X$ on $SM$ we get

$$   \begin{array}{rl}
      \nabla_i f_X  & =  (g_{kl}X^k)_{,i}y^l - g_{kl}X^k \G^l_{ij}y^j, \\
        & = g_{kl}X^k_{,i}y^l + (g_{kj,i} - \G_{k,ij} ) X^k y^j, \\
        & = g_{kl}X^k_{,i}y^l - \frac{1}{2} (g_{kj,i} +  g_{ki,j}  -  g_{ij,k} )  X^k y^j, \\
        & = g_{ka}(X^a_{,i}  +  \G^a_{,li} X^l)y^k, \\
        & = f_{\nabla_i X}.
     \end{array} $$

\qed
Denote the tangent bundle without the zero section by $T'M$. Then $T'M = SM \times \bR^+$ with coordinate
$r$ on the second factor. We denote $E = r \p/ r = y^i \p / \p y^i$ the Euler vector field.
Note that the symplectic form $\o^g$ is homogeneous of degree one: $E \cdot \o^g = \o^g$ and $E \cdot (\o_g)^{-1} = -\o_g $.
Denote the homogeneous functions of degree $k$ by ${\cal{F}}_{k} (T'M)$ and the extension of
the function  $ f \in {\cal F}(SM) $  to ${\cal F} _k$ by  $f_{(k)}$. Furthermore, denote the Poisson
structure by $\{f,g \} = \o_g ^{-1}(df,dg)$. Since $\o_g ^{-1}$ has degree $-1$, we have
$$
\{ {\cal F}_k, {\cal F}_l \} \subset  {\cal F}_{k+l -1}.
$$
We identify ${\cal F}(SM)$ with ${\cal{F}}_{1} (T'M)$, $f \mapsto \tilde{f} = f_1 = f \otimes \bR  $
and define the Legendrian bracket in the space $ {\cal F}(SM)$ by
$$
\{ f, h   \} := \{ \tilde f, \tilde h \}_{SM}.
$$
\qed

\subsection{Smooth  structure in the space of geodesics $L(M) $ }

We first consider the topology of  $L(M)$.

\rm Let $(M,g)$ be a smooth complete Riemannian $n$-dimensional
manifold. By {\it geodesic} we mean an oriented maximally extended
geodesic on $M$ and such a geodesic $\gamma$ has a natural
parameterization by arc-length $\gamma=\gamma(s)$ defined up to a
shift $s\rightarrow s+C$. We denote by $\G$ the canonical geodesic vector
field on the unit sphere bundle $SM$.

 The maximal integral curves of   $\G$ through $(x,v) \in SM$ have the  form
$(\g(s), \g'(s)) = (\exp_x(sv), d/ds \; \exp_x(sv)) $, where $\g(s)=
\exp_x(sv)$ is the maximal geodesic defined by $(x,v) \in SM$.
The set $L(M)$ is identified with the set of orbits of the flow
generated by $\G$, i.e. the maximal integral curves of $\G$. Denote by $\pi : SM \ra L(M)= SM/\G$ the natural projection,
equip $L(M)$ with  the weakest topology such that $\pi$ is
continuous, and call the resulting topological space $L(M)$ the
{\it space of geodesics of M}. In general, $L(M)$ is not Hausdorff, but if $M$ is complete and compact, then
$L(M)$ is compact.
\bl The projection $\pi : SM \ra L(M)$ is an open map if $(M,g)$ is
a complete  Riemannian manifold.
\el

\pf  We have to prove that if $U \subset SM$ is an open set  then $\pi(U)$ is open, i.e.
$V:= \pi^{-1}(\pi(U))\subset SM $ is open.
It is  clear since $V = \cup_{t \in \bR}\varphi_t U$ where $\varphi_t = \exp(t\G)$ is
the 1-parameter group of  transformations which
is generated by $\G$.
\qed

To define a smooth structure in $L(M)$ we consider:

\bd i) A function $f$ on an open subset $U \subset L(M)$ is called smooth if its pull back
$\pi^*f$ is a smooth function on $\pi^{-1}U$. We identify the algebra
${\cal F}(U)$ of smooth functions on $U$ with the algebra
${\cal F}(U)^\G$ of $\G$-invariant functions on $\pi^{-1}U$. In particular,
${\cal F}(L(M)) = {\cal F} (SM)^{\G}$. We denote by ${\cal F}_\gamma(L(M))$
the germ of smooth functions at $\gamma$.
\\
ii) A tangent vector $v$ of $L(M) $ at $\g$ is a derivation $v :{\cal F}_\gamma(L(M)) \ra \bR $,
i.e. $v : f \mapsto v \cdot f$ such that $v \cdot (fh) = h(\g) v \cdot f + f(\g) v \cdot (h).$
They give rise to a vector space denoted by $T_\g L(M).$
 \\
iii)  A  vector field on $U \subset L$ is a derivation of the algebra  ${\cal F}(U)$.
We identify the Lie algebra $\cx(U)$ of vector fields with
the Lie algebra $\cx (\pi^{-1}U)^\G$ of $\G$-invariant vector fields on $\p^{-1}U \subset SM$.\\
iv) A k-form $\o \in \O^k(U)$ is a ${\cal F}(U)$- polylinear skew-symmetric map
$$\o : \cx (U) \times ..... \times \cx (U) \ra \cx (U).$$

\ed

\noindent The standard definition of the exterior differential holds on $\O^k.$
Note in particular that  if the  manifold $M$ has a
dense geodesic, then $C^\infty(L)= \bR  $ and there are no
non-trivial vectors  and vector fields on $L$. But if  we  restrict $M$   to a  sufficiently small neighborhoods
$M'$ of  a point, the  algebra $\mathcal F(M')$  will be  non-trivial and we get a non-trivial Lie  algebra of vector fields.

\subsection{Canonical symplectic structure on $L(M)$}

We now give  another, more general, definition  of  tangent vectors in $T_{\g} L(M)$ in terms of Jacobi fields along
$\g$.

\bd
 A  Jacobi tangent  vector  $v \in T_\g L(M)$   at a point   $\g \in L(M)$
is a Jacobi vector  field $Y$ along $\g$  which is normal to $\g$. The space $T^{J}_\g(L(M)) := Jac_{\g}^\perp$ of such
vector fields is  called the Jacobi  tangent vector  space of $L(M)$   at  $\g$.
\ed
Note that $\dim T^{J}_{\g}(L(M)) = 2n-1$. It is useful to give a relation between the tangent  vector space $T_{\g}L(M)$
 and the Jacobi tangent vector space.

 For a tangent vector $Y \in T_mM$ we denote by $Y_{x,y}^v \in V_{x,y}TM $
the vertical lift and by $Y_{x,y}^h \in H_{x,y}TM $ the horizontal lift. We denote by
$\g^S = (\g, \dot \g)$  the natural lift of  a geodesic $\g$  to $SM$.
\bl The horizontal lift $Y^h$ of a Jacobi field  $Y \in Jac^{\perp}_{\g} = T^J_{\g}L(M) $ is a  horizontal $\G$-invariant  vector  field along
$\gamma^S \subset SM$. It defines a  tangent vector
$$\hat Y : \mathcal F(L(M)) \to \mathbb{R},\,\,   f \mapsto Y^h(f). $$
The map $Y \mapsto \hat Y$ is a homomorphism of $T^J_\g(L(M))$ into $T_\g(L(M))$.
\el

Note that Proposition 1.90 in \cite{B} is incorrect: not every Jacobi field is the transverse field
of a geodesic variation, a counterexample being the vertical constant vector field along the
minimal geodesic of  $x_1 + x_2^2 - x_3^2 = 1 $.

We now define the canonical symplectic 2-form $\omega$ on $L(M)$.

\bl
The 2-form $\o =d \theta$ on $SM$ is $\G$-horizontal $(\iota_\G\o = 0 )$ and $\G$-invariant.
\el
\pf
We have $\theta(\G) = 1$ and $\iota_\G \o =0$, therefore
$\G \cdot \o = ( d \iota_\G + \iota_\G d) \o =0$ and $\G \in ker \;\o $.
Therefore $d\theta$ pushes down to a  closed 2-form $\o = \o_L$
on $L(M) = SM / \G.$
\qed

Now  we describe a Poisson structure on $L(M)$, namely

\bl

${\cal F} (L(M)) ={\cal F} (SM)^{\G}$  is a subalgebra of the Lie algebra $({\mathcal F} (SM), \{, \}).$

\el

\pf
It is known that the Hamiltonian field preserves the symplectic form. If $\G \cdot f  = 0 $, then $\G \cdot \tilde f  = 0$
since $\G \cdot  r  =  0$. We then have $\G \cdot \{f,h  \} = \G \cdot \{ \tilde f, \tilde h  \} =
\G \cdot \o_g^{-1}(d \tilde f, d \tilde h  )|SM = 0  . $
We conclude that $L(M)$ has a canonical Poisson structure.
\qed

One can  check that if the form $\o_L$ is non-degenerate, then the  Poisson structure on $L(M)$ is associated with the
symplectic structure $\o_L $, i.e.
$\{f,h  \} = \o_L^{-1}(df, dh), f,h \in {\cal F}(L)={\mathcal F}_1(SM)^{\G}$  if  $\G \cdot f = \G \cdot h =0$,
where
$ \o_L^{-1}(df, dh) := \omega^{-1}(d \tilde f, d\tilde h)$.

\newpage

\section{Pseudo-Riemannian Spaces of Constant Curvature}

\subsection{Spaces of zero curvature}

Let $E = E^{p+1,q}$ be a pseudo-Euclidean vector space of signature
$(p+1,q)$ with basis  $(e_0^+,...,e_p^+,e_1^-,....,e_q^-)$. The
scalar product is given by
$$ g(X,Y) = <x,y> = \sum_0^p x_{+}^i y_{+}^i - \sum_{1}^{q} x_{-}^j y_{-}j.$$

\bd

A vector $v \in E^{p+1,q}$ is said to be {\it timelike} (respectively {\it spacelike}, {\it null}) if
its norm is negative (positive, vanishes). A  straight line is said to be {\it timelike}
(respectively {\it spacelike}, {\it null}) if its tangent vector has that type.

The space of oriented  {\it timelike} (respectively {\it spacelike}) geodesics of $E^{p+1,q}$ is denoted by $L^-(E^{p+1,q})$ (respectively $L^+(E^{p+1,q})$).
Note that by changing the sign of the metric $L^+(E^{p+1,q})=L^-(E^{q,p+1})$.
\ed

 We denote the unit pseudosphere by
$S^{p,q}= \{ x \in E^{p+1,q},\,\,  x^2 = 1 \}$. Here, for $p = 0$ we
assume in addition that $x^0 > 0$  so that $S^{p,q}$ is connected. The induced metric
has signature $(p,q)$ and constant curvature $1$.

We denote by
$SO(E)= SO^0(p+1,q)$ the connected pseudo-orthogonal group which
preserves the scalar product and by  $SO(E)_e= SO^0(p,q)$ the connected
subgroup which preserves the vector $e = e_0^+$.  The group $SO^0(p,q)$
acts transitively on  $S^{p,q}$  and we can identify $S^{p,q}$ with the quotient
$SO^0(p+1,q)/SO^0(p,q)$. Note that the
tangent space $T_e S^{p,q}$ has the orthonormal basis $(e_1^+,...
,e_p^+,e_1^-,....,e_q^-)$

Any non-null oriented  straight line in $E$ can be canonically written
in the form
$$   \ell_{e, v}(t) = \{ v + te  \},$$
where $e$ is  the unit tangent vector (s.t. $e^2 = \pm 1$) and $v$ is a vector orthogonal to $e$.
So we can identify the  space $L^-(E)$  of timelike lines in $E$  with the tangent  bundle
$TS^{p,q}$.  The group $E(p+1,q) = SO^0(p+1,q) \cdot E^{p+1,q} $ of pseudo-Euclidean motions acts in the space
$L^-(E) = T S^{p,q}$ of timelike lines $(e,v): = l_{e,v}$ by

\begin{equation} \label{action}
T_a (e,v) = (e, v + a_{e^\perp}),\,\,   A (e,v) = (Ae,Av),\,\,  a \in E, \, \, A \in SO^0(p+1,q),
\end{equation}
where $a_{e^\perp }  = a -  <a,e>e$.
The group $E(p+1,q)$ also  naturally acts on $S^{p+1,q}$ with the kernel of effectivity $E^{p+1,q}$.

\begin{Prop}  The isometry group $E(p+1,q)$  acts  transitively  on the  space  $L^-(E)$ of straight lines with stabilizer
$SO(n-1) \cdot R^+$ and  this  action commutes  with  the projection $\pi:L^-(E) \to S^{p.q}, (e,v) \to e.$
 \end{Prop}
The proof follows from equation (\ref{action}).

\begin{Prop}
A necessary condition that a subgroup $G \subset E(p+1,q)$ acts transitively on the space
 $L^-(E^{p+1,q})$ of timelike geodesics is that its linear part $LG = G/{G \cap E}$ acts transitively on $S^{p,q}$.
If  the  group $G$ contains the group  $T_E$ of parallel translations  this condition is also sufficient.
\end{Prop}
Proof. The first claim follows from the previous Proposition. Assume now that $G$ contains $E$  and $LG$ acts transitively on $S^{p,q}$. Let $\ell= \ell_{e,v},\,\ell'= \ell_{e',v'}$ be two lines. Using a transformation from $G$ we transform
$\ell_{e',v'}$  into a line $ \ell'' = \ell_{e,v''}$ with the tangent vector $e$ and then using  parallel translation we
transforms $\ell''$ into $\ell$.

\begin{Cor} Let  $E=E^{n}$ be  the Euclidean vector space. Then any connected  subgroup $G$ of the group $E(n)= SO(E) \cdot E$ of Euclidean motions has the form $G = LG \cdot E$ where $LG \subset SO(E)$ is a connected  orthogonal group which acts transitively on $S^{n-1}$, that is $LG$  is one of the groups
$$SO(n),\,U(n/2),\, SU(n/2),\, Sp(1) \cdot Sp(n/4),\, Sp(n/4),$$
$$ G_2, (n= 7),\, Spin(7),\, (n=8),\, Spin(9),\, (n =16).$$
\end{Cor}

We identify the  space $L^-(E^{p+1,q})$  with the homogeneous   space
$$G/H = E(p+1,q)/SO(p,q)\cdot R^*$$
where $H =SO(p,q)\cdot R^+ $ is the  stabilizer of the line $\ell_0 = \ell_{e_0,0}$. Let $E = \mathbb{R}e_0 + W$  be the orthogonal decomposition. We write
the  corresponding reductive decomposition of the homogeneous space $G/H$ as
$$  \gg = \mathfrak{e}(p+1,q)   = \gh + \gm = (\gso (W)) + \mathbb{R}e_0)  + (U + W),$$
where $\gso(E) = \gso(W) + U$ is the reductive decomposition of $\gso(E)$. In matrix notation,
elements  from  the  algebra of Euclidean isometries  $\mathfrak{e}(p+1,q) \subset \ggl(n+1)$ can be written as

$$(A, \lambda e_0, u,w)  =
  \left( \begin{array}{ccc}
     A & u & w \\
      -u^t & 0 & \lambda\\
   0 & 0 & 0
     \end{array}\right), \,  \lambda \in \mathbb{R}, \,  u,w \in \mathbb{R}^{p+q}.
$$
The adjoint action of $(A, \lambda e_0)$ on $\gm = U \oplus W = \{ (u,w)  \}$ is given by
$$ \ad_{(A, \lambda)} (u,w) = (Au, Aw + \lambda u),$$
and the bracket of two elements  of $\rm = U \oplus W$  is given by
$$ [(u,w), (u', w')] = -(u\wedge u', (u\cdot w' - u' \cdot w )e_0) \in \gh,
$$
where $u\wedge u' $ is an element from $\gso(W)$  and dot means the standard scalar product of vectors from $\mathbb{R}^{n-1} = U = W$.
Note that  $L^-(E^{p+1,q} ) = G/H$  is a symmetric manifold since $\gg = \gh + \gm$ is a symmetric decomposition.

The isomorphism $ \ad_{e_0} : U \to W $ allows one to identify $ U $ with $ W $ and the tangent space $\gm = U + W $  with a tensor product
$ \gm = W \otimes {\mathbb{R}}^2$, where $ U = W \otimes f_1 $, $ W = W \otimes f_2 $ and $ f_1, f_2 $ is the standard basis of $ \mathbb{R}^2 $.
The isotropy representation  $\ad_{\gh}$ preserves  the  Grassmann structures $ \gm = W\otimes \mathbb{R}^2 $
and  $\gso(W)$ act on the first factor $W$  and $\lambda e_0$ acts  on the second factor $\mathbb{R}^2$ by the matrix
$$  \lambda \;\ad_{e_0}  =
 \left( \begin{array}{cc}
     0 & 0 \\
      -\lambda & 0
     \end{array}\right).
$$

Moreover, the isotropy action $\ad_{\gh}|_{\gm}$ preserves  the metric $g^W = g|_{W}$
of signature $(p,q)$ in $W $   and the  symplectic structure $\omega_0 = f_1 \wedge f_2$ in $\mathbb{R}^2$.
The tensor product  $ \omega^{\gm} = g^W \otimes \omega_0 $
defines a non-degenerate  $\ad_{\gh}$-invariant 2-form in
 $\gm$ which is extended to an invariant symplectic form $\omega$ in $L^-(E^{p+1,q})$.  The form  $\omega$  is closed since it is  invariant  and
the manifold  $L^-(E^{p+1,q})$  is a symmetric space.
In the case of dimension $n =3$  the  action  $\ad_{\gso(W)}$ on $W$ also preserves  a 2-form $\omega^W$ (which is the volume form of $W$).
Hence we get  an invariant metric $g^{\gm} = \omega^W \otimes \omega_0 $ on $\gm$ which extends  to an invariant  pseudo-Riemannian metric
$g$ of  signature $(2,2) $ on $L^-(E^{p+1, q})$. The quotient $J = g^{-1} \circ \omega$ is  an invariant (hence, integrable)
 complex structure and  the pair $(g,J)$ is an invariant  K\"ahler structure.
Summarizing, we get ({\it cf.} \cite{S2}):

\begin{Th} \label{t:mt1a}
The space $L^-(E^{p+1,q}) = E(p+1,q)/ SO(p,q) \cdot {\mathbb{R}^+}$ is  a symplectic symmetric space with an invariant Grassmann structure defined by  a decomposition $\gm = W \otimes {\mathbb{R}}^2$. Moreover, if  $n = p+1 +q =3$,
it has an invariant K\"ahler structure $(g,J)$ of  neutral signature $(2,2)$.
\end{Th}

\subsection{Spaces of constant non-zero curvature}
We now describe  the  space of oriented timelike and spacelike geodesics  of the pseudo-Riemannian   space  $S^{p,q}$
of constant curvature 1.  Any such geodesic
through $e\in S^{p,q}$ in direction  of  a unit  vector $e_1^{\pm}$  with $(e_1^{\pm})^2 = \pm 1$ is given by
\[
\g^+ = \g_{e_1^+}^+ = cos(s)e + sin(s)e_1^+ \;,\;\;\; \g^- = \g_{e_1^-}^- = ch(s)e + sh(s)e_1^- .
\]
\noindent The subgroup of the stability group $ SO^0 (p,q)$
preserving the spacelike  geodesic $\g^+$ is $SO(p-1,q)$ and the timelike
geodesic   $\g^-$ is   $SO(p,q-1)$. The one-parameter subgroup $SO(2)$ generated by the
element $e \wedge e_1^+ $ preserves $\gamma^+$ and the one-parameter subgroup $SO(1,1)$ generated by  $e \wedge
e_1^-$ preserves $\gamma^-$. Since the group  $SO^0(p+1,q)$ acts
transitively on the space $L^+(S^{p,q})$ of spacelike geodesics
and on  the space $L^-(S^{p,q})$ of  timelike geodesics, we can
represent these spaces as
\[
 L^+(S^{p,q}) =
SO^0(p+1,q)/SO(2) \cdot SO(p-1,q),
\]
\[
L^-(S^{p,q}) =
SO^0(p+1,q)/SO^0(1,1) \cdot SO(p,q-1).
\]
\noindent To get the reductive decomposition associated with
these spaces, fix the orthogonal decomposition
\[
E^{p+1,q} = \bR e \oplus  \bR e^{\pm}_1 \oplus  V^{\pm},
\]
where $ V^+$ is the vector space  of signature $(p-1,q)$ with basis
$(e_2^+,...,e_p^+,e_1^-,...,e_q^-)$ and   $ V^-$ is the vector space
of signature  $(p,q-1)$ with  basis
$(e_1^+,...,e_p^+, e_2^-,...,e_q^-)$.
 Using the metric, we identify the Lie
algebra $\gso(p+1,q)$ of $SO^0(p+1,q)$ with the space of
bivectors $\Lambda^2(E^{p+1,q})$. Then the reductive decomposition
associated with the unit  sphere bundles
$$ S^{\pm}(S^{p,q})= \{ e_1^{\pm} \in T_eS^{\p,q},\, < e_1^{\pm}, e_1^{\pm}> = \pm 1 \}
$$
is given by
$$  \gso(p+1,q) =
\Lambda^2(V^\pm) \oplus (e \wedge V^\pm \oplus e_1^\pm  \wedge V^\pm ) \oplus \bR (e \wedge e_1^{\pm}).
$$
 The bivector $e \wedge e_1^{+}$ (resp.,$e \wedge e_1^{+}$ ) is invariant under the stability subgroup $SO(p-1,q)$
  (resp., $SO(p,q-1)$ )and defines
an invariant vector field $\G$ on $S^{\pm}(S^{p,q})$, which is the geodesic field. It
 is  the velocity field  of  the right action of the subgroup $SO^{\pm}(2) = SO(2), \, SO(1,1)$ of
$SO(p+1,q)$. The  space of geodesics  is the quotient
 $$L^{+}(S^{p,q}) = SO(p+1,q)/SO(p-1,q)\cdot SO^{\pm}(2),\,\, L^{-}(S^{p,q}) = SO(p+1,q)/SO(p,q-1)\cdot SO(1,1) .$$
The corresponding reductive decomposition may be
written as
\[
\gso(p+1,q) = \gh^{\pm} \oplus  \gm ^{\pm} = \bR (e \wedge e_1^{\pm}) \oplus
\Lambda^2(V^\pm) \oplus (e \wedge V^\pm \oplus e_1^\pm  \wedge V^\pm) .
\]
We identify  $\gm^\pm$ with the tangent space
$T_{\gamma^{\pm}}(L^\pm S^{p,q})$. There is also a natural identification with the tensor product
$\gm^\pm = H \otimes V^\pm  $, where $H = {\rm span} (e, e_1^\pm) \simeq \bR^2$
is the 2-dimensional oriented pseudo-Euclidean vector space.
Then the action of the isotropy subalgebra $\gh^\pm$ takes the form:

$$ ad_{e \wedge e_1^{\pm} }: e \otimes x \mapsto
(e \wedge e_1^{\pm})e\otimes x  = -e_1^{\pm} \otimes x, $$

$$ ad_{e \wedge e_1^\pm} : e_1^\pm \otimes x \mapsto
(e \wedge e_1^\pm)e_1^{\pm})\otimes x  = \pm e \otimes  x,$$

$$ ad_{a\wedge b} : e'\otimes x \mapsto
(e' \otimes (a \wedge b)x =  e' \otimes <b,x>a - <a,x> b,$$
for all $a,b,x \in V^\pm$ and $e' \in \bR^2 $.

Note that $L^{\pm}(S^{p,q})$ is identified with the Grassmanian
$Gr_2^{\pm}(\bR^{p+1,q})$ of two-planes of signature $(2, 0)$ or $(1,1)$ and the decomposition
$$ T_{\gamma^{\pm}}(L^\pm S^{p,q}) = H \otimes V^{\pm}, $$ defines  an invariant  Grassmann
structure in  $L^{\pm}(S^{p,q})$.

Denote  by $H^\e$ two-dimensional vector space  with a scalar product $g_H = g_H^{\e}$  of signature
$(2,0)$ for $\e =+$  and $(1,1)$ for $\e =-$, and by $I_H=I_H^{\e}$ the $SO(H^\e)$-invariant
 endomorphism
of $H^\e$ with $I_H^2 = -\e 1$ and by $\o_H = \o_{H^\e}= g_\e \circ I_H $ the invariant volume form.
Let $(V,g_V)$ be a pseudo-Euclidean vector space of dimension $m$. If $m=2$, we denote by $I_V$ the
$SO(V)$-invariant  endomorphism with $I_V^2 =-1$ for  signature $(2,0)$ or $(0,2)$ and with $I_V^2 =1$
for the signature $(1,1)$. Denote also by $\o_V = g_V \circ I_V$ the  volume form of $V$.
\bl
\begin{enumerate}
\item[i)]  Any $SO(H^\e) \times SO(V)$-invariant  endomorphism  of the   space $W^\e = H^e \otimes V$ has the form $A = 1 \otimes A + I_H^\e \otimes B  $ where $A,B \in \ggl(V)^{SO(V)}$ are invariant endomorphisms of $V$.
\item[ii)] Any invariant endomorphism $I$ of $W^\e$ different from $1$ with $I^2 = \pm 1$  is given
(up to a sigh)  by
$I_{\e} = I_H \otimes 1 $ if $m>2$ and by
$I_\e, I'_\e := 1 \otimes I_V, I''_\e = I_H \otimes I_V $ if $m = 2.$
\item[iii)] Any invariant metric on $W^\e$ is proportional to $g := g_H \otimes g_V$ if $m>2$ and is a
 linear combination of the metric $g$ and the neutral metric $g' := \omega_H \otimes \o_V$ otherwise.
\item[iv)]The space of invariant 2-forms has  the basis $\o = \o_H \otimes g_V$ if $m>2$  and
$\o, \o':= g_H \otimes \o_V$  if $m=2$.
\item[v)] The endomorphisms $I,I'$  are skew-symmetric  with respect to  any  invariant metric $h$ on $W^\e$,
hence define a Hermitian or para-Hermitian structure,   and  the endomorphism $I''$ is symmetric
with  respect  to any invariant metric $h$.
\end{enumerate}
\el

Note that the  tensor product of two complex  or two para-complex structures is a para-complex structure
and  the tensor product of a complex and a para-complex structures is a complex structure.

\pf  To prove part i), it is sufficient to write  the endomorphism $A$ in block matrix form with respect to
the decomposition $W = h_1 \otimes V + h_2 \otimes V$,  where $h_1,h_2$ is an orthonormal basis of
$H^+$ or isotropic basis of $H^-$  and  write the  conditions that it is
 $SO(H) \times SO(V)$-invariant. Since the only invariant  endomorphism of $V$ is a scalar if $m>2$
and  is a linear combination of $1, I_V$ if $m=2$, part ii) follows from part i). Parts iii) and iv)  follow
 from the fact that the space of symmetric bilinear forms
$$ S^2(H \otimes V) = S^2(H) \otimes S^2(V) + \Lambda^2H \otimes \Lambda^2(V),$$
and that the space of 2-forms
 $\Lambda^2(H \otimes V)= \Lambda^2(H) \otimes S^2(V) + S^2(H) \otimes \Lambda^2V$.
Now part v) follows from parts iii) and iv).
\qed

Since  the   spaces
\[
L^+(S^{p,q})= SO^0(p+1,q)/SO(2) \cdot SO^0(p-1,q)
\]
\[
L^-(S^{p,q})= SO^0(p+1,q)/SO^0(1,1) \cdot SO^0(p,q-1)
\]
of  spacelike and timelike geodesics
are  symmetric spaces,  any Hermitian pair $(h,I)$  which consists of invariant pseudo-Euclidean
 metric on $\gm^\e $  and  skew-symmetric invariant complex or para-complex structure $I$ (such that
$I^2 =-1$ or $I^2 =1$)  defines an invariant K\"ahler or para-K\"ahler  structure on $L^\pm S^{p,q}$.

We get the following theorem.

\bt \label{t:mt1b}
\begin{enumerate}
\item[i)]
Let $L^+(S^{p,q}) $, (respectively $L^-(S^{p,q})$)
be the space of spacelike (respectively timelike) geodesics in $S^{p,q}$ and
 $ p+q > 3$. Then there exists a unique (up to  scaling)
invariant symplectic structure $\omega$  and a unique (up to sign)
invariant complex structure $I^+ =J$ on $L^+(S^{p,q})$
(respectively para-complex structure $K=I^-$  on  $L^-(S^{p,q})$).
  There exists unique (up to a scaling ) invariant
pseudo-Riemannian metric $g^\e = \o \circ I^\e$  on $L^\e(S^{p,q})$ which is K\"ahler of signature
 $(2(p-1),2q)$ for
$\e =+$ and para-K\"ahler (of neutral signature)  for $\e =-$.
\item[ ii)] If $p+q =3$,  then there  are two linearly independent invariant (parallel and closed)
 2-forms $\o, \o'$ on $L^\e(S^{p,q})$ with values $\o_{\gm} = \o_H \otimes g_V, \, \o'_{\gm} = g_H \otimes \o_V$.
Any invariant metric has the form $h = \lambda g + \mu g'$ where $g'$ is a neutral metric with
value $g'_{\gm } = \o_H \otimes \o_V$.    Any metric $h$ is K\"ahler for $\e=1$ (respectively,
 para-K\"ahler for $\e =-$) with respect to the complex (respectively, para-complex)  structure
$I^\e = I^\e_H \otimes 1 $ with the K\"ahler form $h \circ I^\pm = \lambda \o + \mu \o'$.
Moreover, the endomorphism $I' = 1 \otimes I_V$ of $\gm$ defines  an  invariant parallel
 $h$-skew-symmetric  complex  structure of $L^\e(S^{p,q})$ if $\e =1, (p-1,q)= (2,0)$ or $(0,2)$ or
$\e =-$  and $(p,q-1) = (2,0)$ or $(0,2)$  and skew-symmetric parallel  para-complex structure
if $\e =+, (p-1,q)= (1,1)$ or $\e =-, (p,q-1) = (1,1)$. The K\"ahler or para-K\"ahler structure
$(h, I')$ has  the K\"ahler form $h \circ I' = \lambda \o + \mu \o'$.
\end{enumerate}
\et
 One can easily check  that the form $\omega$ is the canonical  symplectic form of the space of
 geodesics $L^\e(S^{p,q})$.

\newpage
\section{Rank One Symmetric Spaces of Non-Constant Curvature}

In this section we discuss the invariant geometric structures on
the space $L(M)$ of oriented geodesics   of  a rank one Riemannian symmetric space of non-constant
curvature $M = G/K$, that is for  the projective spaces
\[
\bC P^n =SU_{n+1}/ U_n, \qquad \bH P^n = Sp_{n+1}/ Sp_1 \cdot Sp_n,\qquad \bO P^2 = F_4/Spin_9,
\]
and  the dual  hyperbolic spaces
\[
\bC H^n =SU_{1,n}/ U_n,\qquad \bH H^n = Sp_{1,n}/ Sp_1 \cdot Sp_n,\qquad \bO H^2 = F^{non-comp}_4/Spin_9 .
\]

 In all of these cases, the  space of geodesics is a homogeneous manifold
$L(M) = G/H$ where the stability subgroup $H$ is the same for  the compact and  dual non-compact
case, and is given by
$$ H = T^2 \cdot SU_{n-1},\,\, T^1 \cdot Sp_1  \cdot Sp_{n-1}, T^1 \cdot Spin_7. $$
Moreover, in the case of a classical Lie group $G$, the  space $L(M) $ is the adjoint orbit
$ L(M) = \Ad_G I^\e$ of the  element
$ I^\e=h_1^\e={\rm{diag}}(I^\e_2,0,0,...,0) ,
$ where
$$
I_2^\e=\left(
    \begin{array}{cc}
      0 & -\e \\
       1 & 0
   \end{array}\right),
$$
and $\e =1$ in the compact case and $\e = -1$ otherwise.
Main Theorems \ref{MT:2} and \ref{MT:3} describe all invariant structures
(symplectic structures, complex  and para-complex, K\"ahler and para-K\"ahler
structures) on the space of geodesics $L(M)$.

We prove these in three stages: first for the complex and
quaternionic projective spaces, then for their hyperbolic
counterparts and finally for the Cayley projective
plane.

\subsection{Complex projective and hyperbolic spaces}
We now describe the space of real geodesics in complex projective
space
\[
M^1 = {\mathbb C}P^n = SU_{n+1}/U_{n}
\]
and in complex hyperbolic space
\[
M^{-1} = {\mathbb C}H^n = SU_{1,n}/U_{n}.
\]
We set $\mathfrak{g}^1 = \mathfrak{su}_{n+1}$ and
$\mathfrak{g}^{-1} = \mathfrak{su}_{1,n}$. We then choose the
associated reductive decompositions $\mathfrak{g}^{\epsilon}=
j(\mathfrak{u}_n) + \mathfrak{p}^\epsilon $, $\epsilon =
\pm 1$, where

$    j (\gu_{n})=\{\left(
     \begin{array}{cc}
      -tr A & 0 \\
       0    &   A
     \end{array}\right) \left| \quad  A \in\gu_{n} \right. \}
$, \quad $ \gp_{\epsilon}=\{\left(
     \begin{array}{cc}
       0  & -\epsilon X^* \\
       X   &   0
     \end{array}\right) \left| \quad \quad  X \in \mathbb{C}^n  \right. \}
$,
\newline
\noindent and where $X $ is a column vector and $X^*$ denotes the
Hermitian conjugate. We identify $\mathfrak{p}^\epsilon$ with the
tangent space $T_o M^\epsilon, o = e U_n$.
We next describe the stability subalgebra $\gh^\epsilon$ of the
geodesic $\gamma=\rm{exp}(th_1^\epsilon)(o)$, where the
element $h_1^\e$ 
is represented by the matrix
\[
h_1^\e = I^\epsilon=\left(
     \begin{array}{cc}
      I^\epsilon_2 & 0 \\
       0 & 0
     \end{array}\right), \qquad \qquad \rm{and} \qquad\qquad I^\epsilon_2=\left(
     \begin{array}{cc}
      0 & -\epsilon \\
       1 & 0
     \end{array}\right).
\]
We have that
$\gh^\e = Z_{\gg}(h^\epsilon_1)={\mathbb{R}}h^\epsilon_1 + Z_{\gu_{n}}(h^\epsilon_1)$,
where $Z_{\gu_{n}}(h^\epsilon_1)$ is the centraliser of $h^\epsilon_1$
in $\gu_{n}$. We now describe the reductive decomposition
$\gg^\epsilon =\gh^\epsilon + \gl^\epsilon$.
We have
\[
\gh^\epsilon =\{  A = \left(\begin{array}{ccc}
    i\alpha  & -\epsilon \beta & 0 \\
     \beta    &     i\alpha &     0 \\
    0 & 0 & A_{n-1}
\end{array}\right) | \quad \a, \b \in{\mathbb{R}},
\quad A_{n-1}\in\gu_{n-1}, tr A_{n-1} + 2i \a = 0 \},
\]
and the complimentary subspace is
\[
\gl^\e=\{ X = (x_1,x_2,X_1,X_2)= \left(\begin{array}{ccc}
    i x_1 & \e ix_2 & -\e X_1^* \\
    i x_2 & -i x_1 & -X_2^* \\
    X_1 & X_2 & 0
\end{array}\right) | \quad x_1,x_2\in {\mathbb{R}}, \quad X_1,X_2\in {\mathbb{C}}^{n-1}\}.
\]
We may write $$
\gh^\e = \bR h_0 + \bR h^\e_1 + \mathfrak{su}_{n-1},
$$
where

\[
h_0 =              \left(\begin{array}{ccc}
                    i & 0  & 0 \\
                    0 & i & 0 \\
                      0 & 0 & \frac{-2i}{n-1} {\mathrm Id}_{n-1}
                   \end{array}\right),
\quad
h^\e_1 = I^\e = \left(\begin{array}{ccc}
    0 & -\e  & 0 \\
    1 & 0 & 0 \\
    0 & 0 & 0
\end{array}\right).
\]
Similarly,
$$
\gl^\e = V_0^\e + V_+^\e  + V_-^\e,
$$
where
$$
V_0^\e = \{(x_1,x_2,0,0)\} = \bR E_1 + \bR E_2^\e,\;\; V_\pm^\e  = \{X_\pm = (0,0,X,\pm X), X \in \bC^{n-1} \},
$$

\[
E_1 = (1,0,0,0)=  \left(\begin{array}{ccc}
                    i & 0  & 0 \\
                    0 & -i & 0 \\
                      0 & 0 & 0
                   \end{array}\right), \quad
E_2^\e = (0,1,0,0)=  \left(\begin{array}{ccc}
                    0 & \e i  & 0 \\
                    i & 0 & 0 \\
                      0 & 0 & 0
                   \end{array}\right).
\]

We  denote the canonical Hermitian form in the
space of vector columns $\bC^{n-1}$ by $\eta(X,Y) = X^*Y$ . Then
$$ g(X,Y) = Re \; \eta (X,Y) = \frac{1}{2}(X^*Y + Y^*X ),
$$
$$
\rho (X,Y)= Im \; \eta (X,Y) = \frac{1}{2i}(X^*Y - Y^*X ).
$$
For any $X \in \bC^{n-1} $ we set $X_\pm = (0,0,X,\pm X).$

\bl

We have the following  commutator relations:
 \be
 [E_1, E_2^\e] = 2 h_1^\e,
 \ee
 \be
 [E_1, (0,0,X_1,X_2)] = (0,0,-iX_1, iX_2),
 \ee
\be
 [E_2^\e, (0,0,X_1,X_2)] = (0,0,-iX_2, -\e iX_1).
 \ee
\noindent The isotropy action of $\gh^\e$ on $\gl^\e$ is given by
\[
{\rm ad}_{h_0}(x_1,x_2,X_1,X_2)=(0, 0 ,-i X_1,- iX_2),
\]
\[
{\rm ad}_{h_1^\e}(x_1,x_2,X_1,X_2)=(-2\e x_2,2 x_1,-X_2,\e X_1),
\]
\[
{\rm ad}_{A_{n-1}}(x_1,x_2,X_1,X_2)=(0,0,A_{n-1}X_1,A_{n-1}X_2).
\]
Moreover
 \be
 [X_\pm^1 ,Y_\pm^1  ] = 2\rho(X,Y)(-h_0 \mp E_2^1 )\quad \rm{mod} \quad\gsu_{n-1},
 \ee
\be
 [X_\pm^{-1} ,Y_\pm^{-1}  ] = 2\rho(X,Y)(E_1  \pm E_2^{-1} )\quad \rm{mod} \quad\gsu_{n-1},
 \ee
\be
[X_+^1 ,Y_-^{-1}] = -2\rho(X,Y)E_1 + 2 g(X,Y)h_1 \quad \rm{mod} \quad\gsu_{n-1},
 \ee
\be
[X_+^{-1} ,Y_-^{-1}] = 2\rho(X,Y)h_0 + 2 g(X,Y)h_{-1} \quad \rm{mod} \quad\gsu_{n-1},
 \ee
for all
$
X_\pm^\e, Y_\pm^\e \in  V_\pm^\e
.$
\el

\bp The $\mathfrak{h}$-module $\mathfrak{l^\e}$ has the
following decomposition into irreducible components.

a) For $\e = 1$ \\
$$
\begin{array}{cccc}
  \mathfrak{l}_+^1  = &\quad V_0^1 &\quad  + \quad V_+^1 &\quad + \quad V_-^1 \\
 \mbox{ad}_{h_0}: &0  & -i \Id & \quad -i \Id  \\
  \mbox{ad}_{h_1^\e}: & 2J_0 & i \Id &  \quad -i \Id \\
  A_{n-1}: & 0 & A_{n-1} & \quad A_{n-1}.
\end{array}
$$
b) For $\e = -1$ \\
$$
\begin{array}{ccccc}
  \mathfrak{l}_-^{-1}  =  & \quad \mathbb{R}E_+ &\quad +
  \quad \mathbb{R}E_-  &\quad +   \quad V_+^{-1} \quad + & \quad  V_-^{-1} \\
 \mbox{ad}_{h_0}: &0&0  & -i \Id & -i \Id  \\
  \mbox{ad}_{h_1^\e}: & 2&-2&   \Id &  -\Id \\
  A_{n-1}: & 0&0 & A_{n-1} & A_{n-1},
\end{array}
$$
where $$E_\pm =(1, \pm 1, 0,0) = E_1 \pm E_2^{-1}.$$
\ep

 With this notation in case b), the commutation relations read as follows :
\be\label{1}
[E_{\pm}, V_{\pm}] =0,\; [E_+, E_-] = - 4h_1^{-1},
\ee
\be\label{2} [E_+, X_-] = -2i X_+ ,
\;  [E_-, X_+] = -2i X_- \; \ee
\be\label{3}
[X_\pm , Y_\pm ] = 2 \rho(X,Y)E_\pm, \;\;\;\; [X_+, Y_-] =
2\rho(X,Y)h_0 + 2g(X,Y)h_1 \quad
\rm{mod}\;\mathfrak{su}_{n-1}, \ee
$\mathrm{for} \;\; X_\pm \in V_\pm^{-1}$.
\noindent
Recall that any invariant 2-form on $L(M)$ is generated by an
$\ad_\gh$-invariant two form $\o $ on the tangent space
 $\gl = T_0L(M)$. Any such form may be represented as  $ \o = d( B \circ h )$,
where  $h \in Z(\gh)$   is a central element and $B$ is the Killing form.

  \bt \label{t:mt2a}

\begin{enumerate}
\item [i)]  The only invariant almost complex structures on
    $L({\mathbb{C}}P^n)$ are defined by
\[
J_{\e_0\e_1\e_2}=\e_0\frac{1}{2}{\rm ad}_{h_1^1}|_{V^0}\oplus
\e_1{\rm ad}_{h_1^1}|_{V^+}\oplus \e_2{\rm ad}_{h_1^1}|_{V^-}
= \e_0 J^{V_0} \oplus \e_1 J^{V_+} \oplus \e_2 J^{V_-},
\]
where $\e_k=\pm1$. The integrable ones among these are  (up to a sign) $J_\pm = (\pm J^{V_0}) \oplus J^{V_{-}} \oplus J^{V_{+}}.$
 \item [ii)] Any  $SU(n+1)$-invariant closed 2-form on
 $L({\mathbb{C}}P^n)$ is a linear combination $\o = \l_0\o_0 + \l_1 \o_1 $ of the invariant differential forms
defined by
 $\o_0  = d (B \circ h_0)|\gl^1 $ and $\o_1= d (B \circ h_1^1)|\gl^1$ where $B$ is the Killing form.

Moreover
$$ker \; \o_0 = V_0,\quad
 \o_0 (X_\pm, Y_\pm)= 2 \rho(X,Y),\quad \o_0 (V_+^1,V_-^1 ) = 0,
 $$
 $$
\o_1 (E_1,E_2) = -2, \quad \o_1 (V_\pm,V_\pm) = 0, \quad \o_1
(X_+, Y_-)= -\o_1 (Y_-, X_+) = 2g(X,Y).
 $$
 Here $X_\pm = (0,0,X,\pm X)$, $Y_\pm = (0,0,Y,\pm Y) \in \mathfrak{l^1}$
 and  $$g(X,Y)= \mathrm{Re}\; \eta(X,Y), \o(X,Y) = \mathrm{Im}\; \eta(X,Y)$$
 are the real and imaginary parts of the standard Hermitian form
 $\eta(X,Y) = X^* Y$ on $\mathbb{C}^{n-1}.$
 \item [iii)] The canonical
symplectic structure on the space oriented geodesics is $\o_1$.

\end{enumerate}
 \et
\bc Up to scaling, any invariant symplectic form on
$L({\mathbb{C}}P^n)$ may written as
$$
\o^{t} =   \o_1 + t \o_0, \;\; t \in \mathbb{R}
$$
It is compatible with any invariant complex structure $J_\pm$.
That is, the pair $(\o^{t}, \pm J_\pm)$ is a K\"ahler structure on
$L$.
 \ec
\noindent
Proof of the Theorem: The description in part i) follows directly from the previous
Proposition. For integrability, we calculate the Niejenhuis bracket
$$
N_J(X,Y) = [JX, JY] - J[X,JY] - J[JX,Y] - [X,Y],
$$
for
$X, Y \in \gl = \mathbb{R}E_1 + \mathbb{R}E_2^\e + V_+ + V_-$. For
example, for $X_+ \in V_+,Y_- \in V_- $, $J = J_{\e_01 \epsilon}$ we
calculate :
$$
\begin{array}{cc}
N_J(X_+,Y_-) &= [JX_+, JY_-] - J[X_+,JY_-] - J[JX_+,Y_-] - [X_+,Y_-] \\
& = [i X_+ ,\e i Y_- ] - J [ X_+ , \e i Y_- ] - J[i X_+ , Y_- ] - [X_+,Y_-] \\
&= (-\e  - 1)  [X_+,Y_-] +  (-\e  - 1)i J  [X_+,Y_-]. \\
\end{array}
$$
Hence  $N_J(X_+,Y_-) =0$ iff $J = J_{\e_0 1 1}$ i.e. $\e =1$.
Similarly we calculate $N_J(X_+,E_1)$ :
$$
\begin{array}{cc}
N_J(X_+,E_1) &= [JX_+, JE_1] - J[X_+,JE_1] - J[JX_+,E_1] - [X_+,E_1] \\
&= [i X_+ ,  E_2 ] - J [ X_+ ,  E_2 ] - J[-i X_+ , E_1 ] - [X_+, E_1] \\
&=  (i -J)[X_+,E_2] +  (iJ  - 1) [X_+,E_1] \\
&= (-i + J)i X_+     (iJ  - 1) X_+
\end{array}
$$
This always vanishes. The proof of part ii) follows from the fact that the center
$Z(\gh^1)$ has basis $h_0 , h_1^1.$ To verify part iii) we consider
$h_1^1 \in \gp^1 = T_0\bC P^n.$ The stability subalgebra of this element is
$\gh_{h_1^1}= \bR h_0 + \mathfrak{ su}_{n-1} \subset \gh^1$.
The sphere bundle $S\bC P^n$ is identified with $SU_{n+1}/T^1 \cdot SU_n$ with the
reductive decomposition $\mathfrak{su}_{n+1} = (\bR h_0 + \mathfrak{su}_{n-1}) + (\bR h_1^1  + \gl^1 ).$
The geodesic vector field $\G$ on $ S\bC P^n$ is the invariant vector field
generated by the element $h_1^1$ and the contact form $\theta$ is the invariant form
associated the one-form $\theta_0 = cB \circ h_1^1 $. This shows that the canonical
form coincides with $\o^1$ (up to scaling).

\qed

\bt \label{t:mt3a}
\begin{enumerate} \item [i)] There is no invariant almost
complex structure on the space $L(\mathbb{C}H^n)$. There exist two
(up to sign) almost para-complex structures $K^\pm$ with
$(\pm 1)$-eigenspace decompositions given by
$$
 K^+ : \gl^{-1} = \gl_+ + \gl_- = (\mathbb{R}E_+ + V_+) + (\mathbb{R}E_- + V_-), $$$$
 K^- : \gl^{-1} = \gl_+ + \gl_- = (\mathbb{R}E_+ + V_-) + (\mathbb{R}E_- + V_+).
$$
Only $K^+$ is integrable.
\item [ii)] Any closed invariant two-form is a linear combination
     of the form defined by

$$\o_0 = d (B \circ h_0)| \gl^{-1},\;\; \o_1 = d (B \circ h_1^{-1})| \gl^{-1},  $$
where $B$ is the Killing form.

Moreover, we have
$$\rm{ker}\,\o_0 = V_0,\,\, \o_0(V_\pm, V_\pm) = 0 ,\;\;  \o_0(X_+, Y_-)  =
2 \rho(X,Y),$$
$$\o_1(E_+, E_-) = 4,\; \o_1(X_+, Y_-) =
-\o_1 (Y_-, X_+) = 2g(X,Y),$$
where $X_+ = (0,0,X,X), \; Y_- = (0,0,Y,-Y).$
\end{enumerate}
\et
\noindent
Proof: An  invariant almost complex structure on $L(\mathbb{C}H^n)$
 preserves the one-dimensional $ad_{h_1^{-1}}$ - eigenspaces
 $\mathbb{R}E_\pm$, which is impossible. Since $\gl_\pm = \mathbb{R}E_\pm + V_\pm$
 are subalgebras, the endomorphisms $K$ of $\gl^{-1}$ with $K|\gl_\pm =
 \pm \mathrm{Id}$ define an invariant para-complex structure. It is unique up
 to sign since $[V_\pm , V_\pm] = \mathbb{R}E_\pm$. This proves i).
The first claim of ii) follows from the remark that
  $h_0, h_1^{-1} $ form a basis of the center $Z(\gh^{-1})$. The explicit
  formulas for $\o_0 , \o_1$ follow from equations (\ref{1}), (\ref{2}) and (\ref{3}).
\qed

 \bc Any invariant symplectic form on $L\mathbb{C}H^n$ may
be written as
$$ \o^t = \o_1 + t \o_0.
$$
They are compatible with the para-complex structures $K^{\pm}$, i.e.
$(\o^t,K^{\pm})$ is a para-K\"ahler structure and, in particular, $g^t = \o^t
\circ K^{\pm}$ is a para-K\"ahler metric. \ec

\subsection{Quaternionic projective and hyperbolic spaces}
Consider now the spaces $M^+ = \mathbb{H}P^n = Sp_{n+1}/Sp_1 \cdot
Sp_n$ and $M^- = \mathbb{H}H^n = Sp_{1,n}/Sp_1 \cdot Sp_n$. The
reductive decomposition $\gg^\e = (sp_1 + sp_n) + \gp^\e$ associated
to the homogeneous space  $M^\e$ may  be written as
\[
sp_1 + sp_n = \{\left(
     \begin{array}{cc}
      a & 0 \\
       0 & A_{n}
     \end{array}\right) \left| \; a\in{\rm{Im}}\bH=\gsp_1 , \;\;
 A_{n}\in\gsp_{n} \right. \} , \;
 \gp^\e =\{\left(
     \begin{array}{cc}
      0 & -\e X^* \\
       X & 0
     \end{array}\right) \left | \quad X\in{\mathbb{H}}^n \right\}.
\]
We next describe the stability subalgebra $\gh^\epsilon$ of the
geodesic $\gamma=\rm{exp}(th_1^\epsilon)(o)$, $o = e(Sp_1 \cdot
Sp_n) \in M^\e$, which is the orbit of the one-parameter group
$\rm{exp}(th_1^\e)$,  where $h_1^\e = \mathrm{diag}(I_2^\e,0)$.
We have $\gh=Z_{\gg}(h_1^\epsilon)=
{\mathbb{R}}h_1^\epsilon+Z_{(sp_1+ sp_n)}(h_1^\epsilon)$,
 where $Z_{(sp_1+ sp_n)}(h_1^\epsilon)$ is the
centraliser of $h_1^\epsilon$. 

The reductive decomposition associated to $L(M)^\e$ may be written
as
\[
\gh^\epsilon =\{  A = \left(\begin{array}{ccc}
    a  & -\epsilon \a & 0 \\
     \a    &     a &     0 \\
    0 & 0 & A_{n-1}
\end{array}\right) = ah_0 + \a h_1^\e + A_{n-1}
, a \in sp_1 ,\; \a \in {\mathbb{R}},  A_{n-1} \in sp_{n-1} \}, \]
where
$
h_0 = \rm{diag}(1,1,0), h_1^\e = \rm{diag}(I^\e_2,0)$,
 and the complimentary subspace is
\[
\gl^\e=\{ X = (x_1,x_2,X_1,X_2)= \left(\begin{array}{ccc}
     x_1 & \e x_2 & -\e X_1^* \\
     x_2 & - x_1 & -X_2^* \\
    X_1 & X_2 & 0
\end{array}\right) | \quad x_1,x_2\in sp_1, \quad X_1,X_2\in {\mathbb{H}}^{n-1}\}.
\]
We set
\[
E_1 = \rm{diag}(1,-1,0)\;\;,\;\; E_2^\e = h_1^\e.
\]
Then
\[
\gl^\e = \mathfrak{sp}_1 E_1 + \mathfrak{sp}_1 E_2^\e  + \{(0,0,X_1, X_2) , X_i \in \bH^{n-1}\}.
\]
The isotropy action of $\gh^\e $ on $\gl^\e$ is given by
\[
{\rm ad}_{ah_0}(x_1,x_2,X_1,X_2)=([a,x_1],[a,x_2] ,- X_1 a,- X_2a),
\]
\[
{\rm ad}_{h_1^\e}(x_1,x_2,X_1,X_2)=(-2\e x_2,2 x_1,-X_2,\e X_1),
\]
\[
{\rm ad}_{A_{n-1}}(x_1,x_2,X_1,X_2)=(0,0,A_{n-1}X_1,A_{n-1}X_2).
\]
As in the complex case, we introduce the canonical Hermitian form $\eta (X,Y) = X^* Y$ on
$\bH^{n-1}$. We have the following commutator relations:
\[
[x_1E_1, x_2E_2^\e] = -(x_1x_2 + x_2x_1)h_1^\e = -2Re(x_1x_2)h_1^\e,
\]
\[
[x_1E_1, y_1E_1] = [x_1,y_1]E_1 ,\; [x_2E_2^\e, y_2E_2^\e] =\e[x_2,y_2]h_0,
\]
$
\begin{array}{cc}
&[(0,0,X_1,X_2), (0,0,Y_1,Y_2)] =  -(g(X_2,Y_1) +
g(Y_2,X_1))h_1^\e  \\&  - \big(\e \rho (X_1,Y_1) + \rho (X_2,Y_2)\big)h_0 +
 \big(-\e \rho(X_1,Y_1)+ \rho (X_2,Y_2)\big)E_1 \;
 \rm{(mod}\;\;sp_{n-1}),
\end{array}
$ \\
\[
[x_1E_1, (0,0,X_1,X_2)] = (0,0,-X_1x_1, X_2x_1),
\]
\[
[x_2E_2,(0,0,X_1,X_2)] = (0,0,-X_2x_2, -\e X_1x_2).
\]
For $\e = 1$ we have the following decomposition of the $\gh^1$-module
$\gl^1$ into two irreducible  submodules: $  \gl^1 = V_0 + V_1, V_0 = Im
\mathbb{H} + Im\mathbb{H} = \mathbb{R}^6, V_1 = \mathbb{H}^{n-1} +
\mathbb{H}^{n-1}$. Define an $\gh^1$-invariant complex structure
$J^{V_0},J^{V_1} $ by
 \[J^{V_0} = 1/2 ad_{h_1}|V_0: (x_1,x_2) \to (-x_2,x_1) \quad
 J^{V_1} = ad_{h_1}|V_1: (X_1,X_2) \to (-X_2,X_1).
 \]
The commutator relations imply the following Theorem in a similar
way as in the complex case.

\bt\label{t:mt2b}
\begin{enumerate}\item [i)]
There exist (modulo sign) two invariant almost complex structures on
$L(\mathbb{H}P^n)$ defined by $J_\pm = J^{V_0} \oplus \pm J^{V_1}$.
Only one of them, namely $J_+$, is integrable.
\item [ii)]
Up to scaling, there exists a unique invariant symplectic form $\o$
on $L(\mathbb{H}P^n)$ defined by $\o = d( h_1^\e )^*$. More precisely
we have
\[\o ((x_1,x_2,X_1,X_2),(y_1,y_2,Y_1,Y_2)) =\]
 \[
+2 Re (y_1x_2)+2 Re(y_1x_2) + Re h(X_1,Y_2) - Re h(Y_1,X_2). \]
\item[iii)]
The symplectic structure  $\o$  consistent with the complex structure
$J_+$  and gives rise to a K\"ahler structure.
\end{enumerate} \et
For $\e =-1$ we have the decomposition into irreducible $(\ad_{\gh^{-1}})$ -submodules

$$\gl^{-1} = V_2 + V_{-2} + V_1 + V_{-1} $$
\noindent
such that $\ad_{\gh^{-1}}| V_k  =  k \mathrm{Id}$,
$V_2 \simeq V_{-2} \simeq \mathrm{Im}\bH  \simeq \mathfrak{sp}_1$ and
$V_1 \simeq V_{-1} \simeq \mathrm{Im}\bH^{n-1}$ with the standard action of
$\gh^{-1} = \mathfrak{sp}_1 + \mathfrak{sp}_{n-1}$.
We finally define two $\mathrm{ad}_{\gh}$- invariant para-complex structures
$K_\pm$ on $\gh_{-1}$ by
\[
K_+|_{V_1 +V_2} = 1,\;\;  K_+|_{V_{-1} +V_{-2}} = -1, \;\;  K_-|_{V_1} = 1, \;
\]
\[
K_-|_{V_2} = -1,
\;\;  K_-|_{V_2} = -1, \;\;  K_-|_{V_{-1}} = -1, \;\;  K_-|_{V_{-2}} = 1.
\]
We then have the following
\bt\label{t:mt3b}
On $L(\mathbb{H}H^n)$ there exist no invariant almost complex
structures and two (up to sign) unique invariant almost para-complex
structures $K_\pm$, with $K_+$ being integrable and $K_-$ being non-integrable.
\et

\subsection{Cayley projective plane}

Let $M = \mathbb{O}P^2=  F_4/Spin_9$ be the octonian projective
plane and $\gg = \gf_4 = spin_9 + \gp$ the associated reductive
decomposition.
The isotropy group $Spin_9$ acts on the $16$-dimensional  tangent space
$\gp = T_oM$ by the spinor
representation  with  $15$-dimensional spheres as orbits. Let $\g
= \rm{exp(t\,h_1))o}, h_1 \in \gp$ be the geodesic  through the point
$o = e (Spin_{16}) \in M$.

The stability subgroup of $\g$ is $H= SO_2 \cdot
Spin_7$ and the stability subalgebra
$\gh = \mathbb{R}h_1 + Z_{\mathfrak{spin_9}}(h_1) = \mathbb{R}h_1 + {\mathfrak{spin}}_7$.
 We
identify the space of
geodesics in $M=\mathbb{O}P^2$ with $L(M) = F_4/SO_2 \cdot Spin_7$.

Following \cite{GOV} we choose the root system
$R = \{ \pm \e_{i}, \pm \e_i \pm \e_j, 1/2(\pm \e_1 \pm \e_2 \pm
\e_3 \pm \e_4 )\}$
of the complex Lie algebra $\gf_4$ with
respect to a Cartan subalgebra $\mathfrak{a}$ and a system of
simple roots as follows
$$\a_1 = 1/2(\e_1 - \e_2 - \e_3 - \e_4 ),\, \a_2 = \e_4,\,  \a_3 = \e_3 -
\e_4,\,  \a_4 = \e_4 - \e_3. $$
Here, $ \e_i , \;\;i = 1,..., 4$ is an  orthonormal basis of the real space
$\ga_{\bR}  = B^{-1} \circ \mathrm{span}_\bR R $.
We may assume that $ \frac{1}{2}d = -i h_1  \in
\mathfrak{a}^{\bR}$ is the vector dual to the fundamental weight $\pi_1 =
\e_1$. Then the adjoint operator $ad_{d}$ defines a gradation
$$ \mathfrak{f}_4   = \mathfrak{g}_{-2} +  \mathfrak{g}_{-1}+ \mathfrak{g}_{0}+ \mathfrak{g}_{1}
+ \mathfrak{g}_{2},  $$ where $\mathfrak{g}_0 =
Z_{\mathfrak{f}_4}(h_1) =   \mathbb{C}h_1 + \mathfrak{spin}_7^{\mathbb{C}},$ and
$ \mathfrak{spin}_7^{\mathbb{C}}$ has the root system given by
 $$\{ \pm \e_i,\pm \e_i  \pm \e_j , i,j =  2,3,4 \}.$$
  The space $\mathfrak{g}_{\pm 1}$
 is spanned by the root vectors with roots $1/2(\e_1  \pm \e_2  \pm \e_3  \pm \e_4
)$ and
 $\mathfrak{g}_{\pm 2}$ is spanned by the root vectors
with roots $\pm \e_1, \pm( \e_1 \pm \e_i),\; i = 2,3,4.$
 Let $\t$ be the standard
compact involution of $\mathfrak{f}_4$ such that
$\mathfrak{f}_4^\t$ is the compact real form of $\mathfrak{f}_4$.
Then the reductive decomposition associated  with the space of
geodesics can be
written as
 $\mathfrak{f}_4^\t = \gh + \gl = \mathfrak{g}_0^\t +
\big ( (\mathfrak{g}_{-1} +  \mathfrak{g}_1)^\t + (\mathfrak{g}_{-2} +
\mathfrak{g}_2)^\t $ \big ).
 The decomposition
 $$\gl^{\mathbb{C}} = \gl^{10} + \gl^{01} = (\mathfrak{g}_1 +  \mathfrak{g}_2)
 + (\mathfrak{g}_{-1} +  \mathfrak{g}_{-2}) $$
 defines a unique (up to sign) invariant complex structure $J$ on the space
 of geodesics defined by  $J|_{\gl^{10}}= i \Id,\,\,J|_{\gl^{01}}= -i \Id $.
  The 2-form $\o = d(B \circ h_1)$ associated with the central
 element $h_1  \in Z(\gh) = \mathbb{R}h$ defines a unique
(up to scaling) symplectic form compatible with $J$, where $B$ is the Killing form. We get

\bt \label{t:mt2c}
The space $L(\mathbb{O}P^2)= F_4/SO_2\cdot Spin_7  $ admits a
unique (up to  a sign) invariant complex structure $J$,  unique (up to a
scaling) invariant symplectic structure $\o = d B \circ h_1, h_1 = \in
Z(\gh)$ and a unique
invariant K\"ahler structure $(\o , J)$.
\et

\end{document}